\documentclass[a4paper,10pt]{article}
\usepackage{t1enc}
\usepackage{latexsym}
\usepackage{amsfonts}
\usepackage{amssymb}
\usepackage{amsmath}
\usepackage{indentfirst}
\usepackage{amscd}
\usepackage[all]{xy}
\usepackage{multicol}

\newcommand{\Dem}{{\bf Proof: }}

\newtheorem{Def}{Definition}

\newtheorem{Lem}[Def]{Lemma}
\newtheorem{Teo}[Def]{Theorem}

\newcommand{\cqd}{\hfill$\Box$}

\newcommand{\aut}{\mbox{Aut}}
\newcommand{\End}{\mbox{End}}
\title{The Chern-Connes character for pseudodifferential operators on the sphere}
\author{David P. Dias \and Severino T. Melo\footnote{Partially supported by 
CNPq (Processo 304783/2009-9), Brazil.}}
\date{ }

\begin{document}

\maketitle

\begin{abstract}
We compute the Chern-Connes character (a map from the $K$-theory of a 
C$^*$-algebra under the action of a Lie group to the cohomology of its 
Lie algebra) for the $L^2$-norm closure of the algebra of all classical 
zero-order pseudodifferential operators on the sphere under the canonical 
action of ${\rm SO}(3)$. We show that its image is $\mathbb{R}$ if the trace
is the integral of the principal symbol.
\end{abstract}

\vskip0.3cm
\begin{center}
Mathematics Subject Classification:  46L80 (47G30, 58B34). 
\end{center}
\vskip0.3cm

\section*{Introduction}

Years before his seminal work on noncommutative differential geoemetry, A. Connes defined a generalization of the Chern character for noncommutative C$^*$-algebras under the action of Lie groups, as a map from the $K$-theory of the algebra to the cohomology of the Lie algebra of the group. The dense subalgebra of the elements with smooth orbit then plays the role of the smooth functions in the commutative case. This much simpler and more intuitive generalization of the Chern character has not received much attention, perhaps in view of the much more general and powerful algebraic tools developed in and after \cite{AC2}.

In this paper, we compute the Chern-Connes character of \cite{AC1} for the norm closure of the algebra of classical zero-order pseudodifferential operators on the sphere $S^2$ under the canonical action of ${\rm SO}(3)$. The set of classical pseudodifferential operators that generate this algebra is in fact much smaller than the algebra of all bounded operators on $L^2(S^2)$   with smooth ${\rm SO}(3)$-orbit. For example, the algebra of all zero-order pseudodifferential operators in the standard H\"ormander class has been characterized by M. Taylor \cite{MT} as the bounded operators which are smooth under the action of the conformal group 
${\rm SO}_e(3,2)$, which contains ${\rm SO}(3)$ as a proper subgroup. And H\"ormander operators may lie outside the closure of the classical ones \cite{T}. There are also ``smooth'' operators
which do not even have the pseudolocal property \cite{CoMe}.

The $K$-theory of the algebra considered here is quite simple: $K_1$ vanishes and, modulo torsion, $K_0$ is 
generated by the class of the identity. Since the identity is fixed by the action, only the zero-degree cohomology class of the Chern caracter, which is equal to the trace of the projection, could be nonzero. The Wodzicki trace \cite{W}\ vanishes on the identity (in fact on any pseudodifferential projection on a closed manifold) and it is the unique trace, modulo a multiplicative constant, defined on the algebra of pseudodifferential operators of arbitrary order, if the manifold and its cosphere bundle are connected. For the algebra of zero-order operators, however, the integral of the principal symbol defines a trace which, in the case of $S^2$, is ${\rm SO}(3)$-invariant. For this choice of trace, the Chern-Connes character for our algebra has image $\mathbb{R}$.

In the first section of this paper,  we collect definitions and results from \cite{AC1}. 
In Section 2, we compute the $K$-theory of the cosphere bundle of $S^2$, which has also been computed by F. Rochon using different methods \cite{FR}. In Section 3, we state the main result.

\section{ The Chern-Character for a C$^*$-dynamical system} 

A {\em C$^*$-dynamical system} is a triple $(A, G, \alpha)$, where $A$ 
is a unital C$^*$-algebra, $G$ is a Lie group and the homomorphism
$\alpha:G \rightarrow \aut(A)$, 
$g\mapsto\alpha_g$, is continuous in the sense that 
$\alpha(a):G\ni g\mapsto\alpha_g(a)\in A$\  is continuous 
for every $a \in A$. An $a \in A$ is $C^\infty$ if $\alpha(a)$ is $C^\infty$. 
By the Garding Theorem \cite[Theorem V.1]{DA}, 
$A^\infty = \{ a \in A: \textrm{ a is } C^\infty  \} $ is (norm) dense in $A$. 
Moreover, Connes \cite[Appendix]{AC3} shows that $A^\infty$ is invariant under the 
holomorphic functional calculus, and so the inclusion of $A^\infty$ into $A$ induces a 
K-theory isomorphism
$K_i(A^\infty) \simeq K_i(A)$, $i = 0,1$. 

A representation $\delta$ of $\mathfrak{g}$, the Lie algebra of $G$, in 
the derivations of $A^\infty$ is defined by 
$$ \delta_X(a) = \underset{t \rightarrow 0}{\lim} \frac{\alpha_{g_t}(a) - a}{t}, $$
where $g_0' = X$ and $a \in A^\infty$.
The fact that $\delta$ is indeed a Lie-algebra representation is shown, for example, in
\cite[A.4]{DPD2}. There is a canonical extension of $\alpha_g$ to $M_n(A)$ and of 
$\delta$ to $M_n(A^\infty)$.

One uses  $\delta$  to define the complex 
$\Omega = A^\infty \otimes \Lambda \mathfrak{g}^*$ of left invariant 
differential forms on $G$ with coeficients in $A^\infty$, with exterior 
derivative $d$ such that $\delta_X(a) = da(X)$, 
$d(w_1 \wedge w_2) = d(w_1) \wedge w_2 + (-1)^k w_1 \wedge d(w_2)$ and $d^2(w_2) = 0$, 
for all $a \in A^\infty$, 
$X \in \mathfrak{g}$, $w_1 \in \Omega^k$ and $w_2 \in \Omega$.

A {\em connection} on a finitely generated module $M^\infty$ over $A^\infty$ 
is a linear map 
$\nabla: M^\infty \rightarrow M^\infty \otimes \mathfrak{g}^*$ such that $ \nabla_X(\xi a) = \nabla_X(\xi)a + \xi \delta_X(a),$
for $\xi \in M^\infty$, $X \in \mathfrak{g}$ and $a \in A^\infty$.
Any finitely generated projective module over $A^\infty$  is isomorphic to
$p(A^\infty)^n$, for some $n \in \mathbb{N}$ and some idempotent 
$p\in M_n(A^\infty)$. We can suppose that $p$ is self-adjoint \cite[4.6.2]{BB}.
The {\em grassmaniann connection} on 
$p(A^\infty)^n$  is defined by $\nabla^0_X(\xi) = p \delta_X(\xi)$, 
$\xi \in M^\infty$, $X \in \mathfrak{g}$. 
Associated to a connection $\nabla$ on $M^\infty\simeq p(A^\infty)^n$, 
the {\em curvature}  $\Theta\in\End_{A^\infty}(M^\infty) \otimes \Lambda^2(\mathfrak{g}^*)$ is  defined by 
$$
\Theta(X,Y) = \nabla_X \nabla_Y - \nabla_Y \nabla_X -\nabla_{[X,Y]},\  
X,Y\in\mathfrak{g}.
$$
Identifying $\End_{A^\infty}(M^\infty)$ with 
$pM_n(A^\infty)p \subset M_n(A^\infty)$ by the 
isomorphism $f \mapsto pf(p)p$ and using the fact that $\delta$ is a 
representation, one can check that the curvature associated to 
the grassmaniann connection $\nabla^0$ is the $2$-form 
$\Theta_0 = p dp \wedge dp \in \Omega^2$.

A linear map $\tau: A \rightarrow \mathbb{C}$ is a {\em finite G-invariant trace} 
if it vanishes on commutators, it is positive $($and then 
$\tau(a^*) = \overline{\tau(a)}$ for every $a\in A$),  and if
$\tau(\alpha_g(a)) = \tau(a)$ for every $a\in A$ and $g\in G$.
Given a finite G-invariant trace $\tau$ on $A$, for each positive integer 
$k$, there exists a unique k-linear map 
$\tau_k : \Omega \times \ldots \times \Omega \rightarrow \Lambda \mathfrak{g}^*$
such that 
\[
\tau_k (a_1 \otimes w_1, \ldots , a_k \otimes w_k) = \tau(a_1 \ldots a_k) 
w_1 \wedge \ldots \wedge w_k, \ \ a_i\in \mathfrak{g}^*
\]
We also denote $\tau_0=\tau$.

The differential form 
$\tau_k(\Theta, \Theta, \ldots, \Theta) \in \Lambda^{2k}\mathfrak{g}^*$ 
is closed and its cohomology class does not depend on the choice of connection 
on $M^\infty$ \cite[Proposition 5]{AC1}. It then follows that it also does not 
depend on the choice of a representative in an isomorphism class of finitely 
generated projective modules. Using that $K_0(A^\infty)\simeq K_0(A)$, 
one can then define the {\em Chern-Connes character} as the homomorphism 
$Ch_\tau: K_0(A) \rightarrow H_{\mathbb{R}}^{\mbox{{\scriptsize even}}}(G)$, 
$$ Ch_\tau ([p]_0) = \left[ \underset{k=0}{\overset{\infty}{\sum}} \frac{1}{(2 \pi i)^k k!} \tau_k(\Theta, \ldots, \Theta) \right]=$$
\begin{equation*}\label{chernconnes}
\begin{split} 
= \left[ \underset{k=0}{\overset{\infty}{\sum}} \frac{1}{(2 \pi i)^k k!} \tau_k(\Theta_0, \ldots, \Theta_0) \right] 
= \left[ \underset{k=0}{\overset{\infty}{\sum}} \frac{1}{(2 \pi i)^k k!} \tau_1(p(dp \wedge dp)^k)\right].
\end{split}
\end{equation*}
To show that $Ch_\tau ([p]_0)$ is indeed real, observe that, for 
$da \wedge db \in A^\infty \otimes \Omega^2$ and $X, Y \in \mathfrak{g}$, we have  
$
(da \wedge db)^*(X,Y) = (\delta_X(a)\delta_Y(b)-\delta_Y(a)\delta_X(b))^*
= - (\delta_X(b^*)\delta_Y(a^*)-\delta_Y(b^*)\delta_X(a^*))= - (db^* \wedge da^*)(X,Y).
$
Hence 
$\overline{\tau_{1}(p (dp \wedge dp)^k)} = \tau_{1}((p(dp\wedge dp)^k)^*) =  (-1)^k \tau_{1} (p(dp \wedge dp)^k),$ and 
so $i^k\tau_{1}(p (dp \wedge dp)^k)$ is real.

Given a C$^*$-dynamical system $(A, G, \alpha)$, another C$^*$-dynamical system 
$(A \otimes C(S^1), G \times S^1, \alpha')$ is canonically defined, with action 
$\alpha'_{(g,h)} (a \otimes f) = \alpha_g(a) \otimes f_h$, 
where $f_h(x) = f(x-h)$.
Using the split exact sequence 
\begin{equation}\label{split}
0 \rightarrow SA \rightarrow C(S^1, A) \rightarrow A \rightarrow 0
\end{equation}
where the map onto $A$ is evaluation at $1$ and
$SA$ denotes the suspension of $A$,
we obtain 
\[
K_0(C(S^1)\otimes A)\simeq K_0(C(S^1,A)) \simeq K_0(A) \oplus K_0(SA) \simeq K_0(A) \oplus K_1(A).
\]
By \cite[3B.3]{H1}, $H^k(G \times S^1) \simeq H^k(G) \otimes H^{k-1}(G)$, so the map  
$Ch_\tau: K_0(C(S^1) \otimes A) \rightarrow 
H_\mathbb{R}^{\mbox{{\scriptsize even}}}(G \times S^1)$ can be 
understood as
$ch_\tau : K_0(A) \oplus K_1(A) \rightarrow H_\mathbb{R}^*(G)$.

\section{ K-theory computation}


For any closed manifold, if $\mathcal{A}$ denotes the closure of the algebra of zero-order classical pseudodifferential operators and $\mathcal{K}$ denotes the ideal of compact operators, then the principal symbol induces an isomorphism between the quotient 
$\mathcal{A}/\mathcal{K}$ and the algebra of continuous functions on the cosphere bundle. This fact follows from the classical estimate for the norm, modulo compacts, of pseudodifferential operators \cite[Theorem A.4]{KN}.  In the case considered here, where the manifold is the sphere $S^2$, we then get the exact sequence
\begin{equation} \label{seqS2}
 0 \longrightarrow \mathcal{K} \longrightarrow \mathcal{A} \overset{\sigma}{\longrightarrow} C(S^*S^2)\longrightarrow 0,
\end{equation}
with $\mathcal{A}$ denoting the $L^2(S^2)$-norm closure of the algebra $\Psi_{cl}^0 (S^2)$ of all zero order classical 
pseudo\-differential operators on $S^2$, and $S^*S^2$ denoting the subset of the cotangent
bundle $T^*S^2$ consisting of covectors of norm 1, for some choice of 
Riemannian metric on $S^2$.

Before we use  (\ref{seqS2}) to compute the K-theory groups of $\mathcal{A}$, we need to compute the K-theory of
$C(SS^2) \simeq C(S^*S^2)$. For that, we will use  the following version of the Mayer-Vietoris sequence 
(see \cite[7.2.1]{SM} or \cite[21.5.1]{BB}).

\begin{Teo} \label{MVH}
Given the unital C$^*$-algebra comutative diagram 
$$\xymatrix{ A \ar[d]_{p_2} \ar[r]^{p_1}& B_1 \ar[d]^{\pi_1}\\
B_2 \ar[r]_{\pi_2}& B} $$
where $A$ is the bundle product of $B_1$ and $B_2$ on $B$, that is,  
$$A = \{ (b_1, b_2) : \pi_1(b_1) = \pi_2(b_2) \} \subseteq B_1 \oplus B_2 \ ,$$  
with $\pi_1$ and $\pi_2$  surjective $*$-homomorphisms and $p_k$, $k=1 \textrm{ or } 2$, restrictions of the projections  $i_k : B_1 \oplus B_2 \rightarrow B_k$, then we have the cyclic exact sequence with six terms
$$ \xymatrix{
K_0(A)  \ar[rr]^-{(p_{1_*},p_{2_*})}  & \ \ &  K_0(B_1) \oplus K_0 (B_2) \ \ \ar[rr]^-{\pi_{2_*} - \pi_{1_*}}     & \ \ &   K_0(B) \ar[dd]     \\
& \\
K_1(B) \ar[uu]^{}  & \ \ &  K_1(B_1) \oplus K_1(B_2) \ar[ll]_-{\pi_{2_*} - \pi_{1_*}}  & \ \  & K_1(A) \ar[ll]_-{(p_{1_*},p_{2_*})}                     } .$$

\end{Teo}


In the proof of the following lemma, we use the charts induced on $SS^2$ by the stereographic projections to identify $C(SS^2)$ with a bundle product of two copies of $C(D\times S^1)$, where $D=\{x\in\mathbb{R}^2;\,|x|\leq 1\}$ and $S^1$ is the unit circle. We have denoted 
by $SS^2$ the subset of unit tangent vectors over $S^2$, for some choice of metric. The algebras $C(S^*S^2)$ and $C(SS^2)$ are naturally isomorphic.

\begin{Lem} \label{diagrama} Let $\pi_1:C(D\times S^1)\to C(S^1\times S^1)$ denote the restriction map,
$\pi_1\varphi= \varphi_{ \mid_{S^1 \times S^1}}$, and let $\pi_2:C(D\times S^1)\to C(S^1\times S^1)$
be defined by $(\pi_2\varphi)(z,w) = \varphi (z, -z^2 \bar{w})$, $(z,w)\in S^1 \times S^1$. Then there are 
C$^*$-algebra homomorphisms $p_i:C(SS^2)\to C(D\times S^1)$, $i=1,2$, such that 
$(p_1,p_2):C(SS^2)\to C(D\times S^1)\oplus C(D\times S^1)$ is injective and its image is equal to
$$\{ (\varphi_1,\varphi_2) \in C(D \times S^1) \oplus C(D \times S^1) ;\,\pi_1(\varphi_1)=\pi_2(\varphi_2)\}.$$
\end{Lem}
\Dem 
Define $\chi_N : S^2\setminus\{ (0,0,1) \} \rightarrow \mathbb{R}^2$ and $\chi_S : S^2\setminus\{ (0,0,-1) \} \rightarrow \mathbb{R}^2$ by
\[
\chi_N(x_1, x_2, x_3)=\left( \frac{x_1}{1-x_3} , \frac{x_2}{1-x_3} \right) = (N_1, N_2) 
\]
and
\[
 \chi_S(x_1, x_2, x_3)=\left( \frac{x_1}{1+x_3} , \frac{x_2}{1+x_3} \right) = (S_1 , S_2).
\]
The change of coordinates $K=\chi_S\circ(\chi_N)^{-1}$ is given by $K(N_1,N_2)= \left( \frac{N_1}{N_1^2 + N_2^2} , \frac{N_2}{N_1^2 + N_2^2} \right)$. 

If $T\chi_S$ and $T\chi_N$
denote the charts induced on the tangent bundle by $\chi_S$ and $\chi_N$, then $T=T\chi_S \circ (T\chi_N)^{-1}$ is given by
\[
 T(N_1, N_2, \eta_1, \eta_2) = \left(\frac{N_1}{N_1^2 + N_2^2} , \frac{N_2}{N_1^2 + N_2^2}, \frac{\partial (K_1, K_2)}{\partial (N_1, N_2)} \left( \begin{array}{ll} \eta_1 \\ \eta_2 \end{array} \right) \right),
\]
with
\[
\frac{\partial (K_1, K_2)}{\partial (N_1, N_2)} = \frac{1}{(N_1^2 + N_2^2)^2} \left( \begin{array}{ll}  - N_1^2 + N_2^2 & -2 N_1 N_2 \\ -2 N_1 N_2 &  N_1^2 - N_2^2 \\  \end{array} \right).
\]
Identifying  $\mathbb{R}^2 \times \mathbb{R}^2$ with $\mathbb{C}^2$, and restricting  to $S^1 \times S^1$, we get $T(z,w) = (z, -z^2 \bar{w})$.
The theorem is proved if we define $p_1(f)$ as the restriction of $f\circ (T\chi_N)^{-1}$ to $D\times S^1$ and
$p_2(f)$ as the restriction of  $f\circ (T\chi_S)^{-1}$ to $D\times S^1$  \cqd 

Applying Theorem \ref{MVH} to the algebras and maps of Lemma \ref{diagrama}, we obtain the exact sequence
\begin{equation} \label{CSS2}
 \xymatrix{
K_0(C(SS^2))  \ar[rr]^-{(p_{1_*},p_{2_*})}  & \ \ &  K_0(C(D \times S^1))^2 \ \ \ar[rr]^-{\pi_{2_*} - \pi_{1_*}}     & \ \ &   K_0(C(S^1 \times S^1)) 
 \ar[dd]       \\
& \\
K_1(C(S^1 \times S^1)) \ar[uu]^{}  & \ \ &  K_1(C(D \times S^1))^2 \ar[ll]_-{\pi_{2_*} - \pi_{1_*}}  & \ \  & K_1(C(SS^2)). \ar[ll]_-{(p_{1_*},p_{2_*})}                     } 
\end{equation}

It is well known that $K_0(C(D))=[\mathfrak{1}]_0 \mathbb{Z}$, where $\mathfrak{1}$ denotes
the constant function $z\mapsto 1$, and that $K_1(C(D))=0$; and also that 
$K_0(C(S^1))=[\mathfrak{1}]_0\mathbb{Z}$ and $K_1(C(S^1))=[\mathfrak{z}]_1\mathbb{Z}$, where
$\mathfrak{z}$ denotes the identity function $z\mapsto z$. Using the canonical isomorphisms
$C(D \times S^1) \simeq C(S^1, C(D))$ and $C(S^1 \times S^1) \simeq C(S^1, C(S^1))$ and 
the exact sequence (\ref{split}) with $A$ replaced first by $C(D)$ and then by $C(S^1)$, 
we get
\[
K_0(C(D \times S^1)) = [\mathfrak{1}]_0 \mathbb{Z},   \ \ \ \ 
K_1(C(D \times S^1)) = [\mathfrak{w}]_1 \mathbb{Z},
\]
\[
K_0(C(S^1 \times S^1)) = [\mathfrak{1}]_0 \mathbb{Z} \oplus \theta_{C(S^1)}([\mathfrak{w}]_0) \mathbb{Z}\ \ \ \mbox{and}\ \ \ 
K_1(C(S^1 \times S^1)) = [\mathfrak{z}]_1 \mathbb{Z}\oplus [\mathfrak{w}]_1 \mathbb{Z},
\]
where now $\mathfrak{z}$ and $\mathfrak{w}$ denote, respectively, the maps
$(z,w)\mapsto z$ and $(z,w)\mapsto w$. We use the standard notation $\theta_A$ for the natural isomorphism from $K_ 1(A)\to K_0(SA)$, for a C$^*$-algebra $A$.   

The diagram (\ref{CSS2}) now becomes
$$ \xymatrix{
K_0(C(SS^2))  \ar[rr]^-{(p_{1_*},p_{2_*})_0}  & \ \ &  \mathbb{Z} [\mathfrak{1}]_0 \oplus \mathbb{Z} [\mathfrak{1}]_0 \ \ \ar[rr]^-{(\pi_{2_*} - \pi_{1_*})_0}     & \ \ &   \mathbb{Z} [\mathfrak{1}]_0 \oplus \mathbb{Z} \theta_{C(S^1)}([\mathfrak{w}]_1)  \ar[dd]^{\delta_0}  
     \\
& \\
\mathbb{Z} [\mathfrak{z}]_1 \oplus \mathbb{Z} [\mathfrak{w}]_1  
\ar[uu]^-{\delta_1}
 & \ \ &  \mathbb{Z} [\mathfrak{w}]_1 \oplus \mathbb{Z} [\mathfrak{w}]_1  \ar[ll]_-{(\pi_{2_*} - \pi_{1_*})_1}  & \ \  & K_1(C(SS^2)). \ar[ll]_-{(p_{1_*},p_{2_*})_1}                     } $$

Using the obvious identifications between some of the groups above and $\mathbb{Z}^2$, we have, 
for all $(x,y)\in\mathbb{Z}^2$, 
\begin{equation}
\label{tltc}
(\pi_{2_*} - \pi_{1_*})_0 (x,y) = (y - x,0)\ \ \ \mbox{and}\ \ \ 
(\pi_{2_*} - \pi_{1_*})_1 (x,y) = (-2y, -y-x). 
\end{equation}

Since $(\pi_{2_*} - \pi_{1_*})_1$ is injective, $(p_{1_*},p_{2_*})_1=0$.
And then it follows from the fact that the image of $(\pi_{2_*} - \pi_{1_*})_0$ is $\mathbb{Z}[1]_0$ that 
\[
K_1(C(SS^2))\simeq\mathbb{Z}, \ \ \mbox{with generator}\ \ \delta_0(\theta_{C(S^1)}([\mathfrak{w}]_1)).
\]
It follows from the second equation in (\ref{tltc}) that the quotient of of 
$K_1(C(S^1\times S^1))$ by the kernel of $\delta_1$ is isomorphic to $\mathbb{Z}_2$. 
Since the kernel of $(\pi_{2_*} - \pi_{1_*})_0$ is isomorphic to $\mathbb{Z}$, we 
obtain an exact sequence 
$0\to \mathbb{Z}_2\to K_0( C(S^*S^2))\to \mathbb{Z}\to 0$, and hence
\[
K_0( C(S^*S^2)) \simeq \mathbb{Z} \oplus \mathbb{Z}_2.
\]
It is straighforward to check that the copy of $\mathbb{Z}$ contained in $K_0(C(S^*S^2))$
is generated by $[\mathfrak{1}]_0$. 


It follows from the Fedosov index formula (as in \cite{LM}) that there is a Fredholm operator of index 1 in $\mathcal{A}$. The index map in the standard six-term exact sequence associated
(\ref{seqS2}),
$$ \xymatrix{
\mathbb{Z}  \ar[rr]^-{ 0}  & \ \ &  K_0(\mathcal{A}) \ \ \ar[rr]^-{\sigma_*}     & \ \ &  \mathbb{Z} \oplus \mathbb{Z}_2  \ar[dd]^-{}     \\
& \\
\mathbb{Z} \ar[uu]^-{\delta_1}  & \ \ &  K_1(\mathcal{A})  \ar[ll]_-{\sigma_*}  & \ \  & \ 0, \ar[ll]_-{0}       } $$
is therefore surjective. Then $\sigma_*$ is an isomorphism on $K_0$, and hence
\[
K_0(\mathcal{A})\simeq \mathbb{Z} \oplus \mathbb{Z}_2, 
\]
where the copy of $\mathbb{Z}$ is generated by the class of the identity. The surjective
homomorphism $\delta_1:\mathbb{Z}\to\mathbb{Z}$ is necessarily injective. From this we 
get 
\[
K_1(\mathcal{A})=0.
\]

\section{The Chern-Connes character for $\overline{\Psi_{cl}^0 (S^2)}$}

For each $g \in {\rm SO}(3)$, let $T_g$ denote the unitary operator $T_g u(x) = u(g^{-1}x)$, $u \in L^2(S^2)$. It is well known that, if $A$ is a zero-order classical pseudodifferential operator on $S^2$, then, for each $g\in {\rm SO}(3)$, $T_gAT_g^{-1}$ is also a zero-order pseudodifferential operator. Moreover, the map 
\begin{equation}
\label{action} 
{\rm SO}(3)\ni g\mapsto \alpha_g(A)=T_gAT_g^{-1}
\end{equation}
is smooth with respect to the norm topology of $\mathcal{L}(L^2(S^2))$
(see, for example, \cite{MT}). 
It then follows that $\alpha_g$ extends continuously to an automorphism of 
the algebra $\mathcal{A}$ of the previous section. The smooth algebra
$\mathcal{A}^\infty$ of the C$^*$-dynamical system  $(\mathcal{A},{\rm SO}(3),\alpha)$ then contains
the algebra of classical pseudo\-differential 
operators $\Psi_{cl}^0 (S^2)$.

We saw in the previous section that $K_1(\mathcal{A})=0$ and $
K_0(\mathcal{A})\simeq\mathbb{Z}\oplus\mathbb{Z}_2$. Since any group 
homomorphism from $\mathbb{Z}_2$ to a real vector space must vanish, 
the Chern-Connes character
for this C$^*$-dynamical system  reduces to a map defined on the 
subgroup of $K_0(\mathcal{A})$ isomorphic to $\mathbb{Z}$. We also 
saw in the previous section that this subgroup is generated by 
the class of the identity $[I]_0$. Since $I$ is fixed under the 
action of the group, $dI=0$. Thus the only nonvanishing term in 
the sum (\ref{chernconnes}) is for $k=0$: the constant $\tau(I)$.

The obvious example of an ${\rm SO}(3)$-invariant trace on the 
quotient $\mathcal{A}/\mathcal{K}\simeq C(S^*S^2)$ is the integral 
with respect to the measure induced on $S^*S^2$ by the usual surface
measure on $S^2$. Composing this trace with the quotient projection
$\mathcal{A}\to\mathcal{A}/\mathcal{K}$, we get an ${\rm SO}(3)$-trace on 
$\mathcal{A}$ for which $\tau(I)\neq 0$ ($\tau(I)$ is equal to the 
measure of ${\rm SO}(3)$). This shows that the image of $Ch_\tau$ is isomorphic to 
$\mathbb{R}$.


\vskip0.5cm

\noindent{\footnotesize 
Instituto de Matem\'atica e Estat\'{\i}stica, Universidade de S\~ao Paulo, 
Rua do Mat\~ao 1010,
05508-090 S\~ao Paulo, Brazil.\\
Emails: dpdias@ime.usp.br,\ \ \ toscano@ime.usp.br.}

\end{document}